\documentclass[a4paper,12pt,leqno]{article}

\usepackage{geometry} 

\usepackage{pgfplots}
\usepackage{tikz-cd
}\usepackage{bussproofs}
\usepackage{amssymb,amsmath}
\usepackage{proof}
\usepackage{psvectorian}
\usepackage{comment}
\usepackage{CJKutf8}
\usepackage[utf8]{inputenc}
\usepackage{multicol}
\usepackage{graphicx}
\graphicspath{ {images/} }
\usepackage{qtree}
\usepackage{mathtools}
\usepackage{hyperref}
\usepackage{listings}
\usepackage{indentfirst}
\usepackage{titlesec}
\usepackage{enumerate}

\DeclareUnicodeCharacter{00AC}{$\neg$}

\DeclareUnicodeCharacter{03B5}{$\epsilon$}

\DeclareUnicodeCharacter{2200}{$\forall$}

\DeclareUnicodeCharacter{2203}{$\exists$}

\DeclareUnicodeCharacter{2229}{$\cap$}

\DeclareUnicodeCharacter{222A}{$\cup$}

\DeclareUnicodeCharacter{2282}{$\subset$}

\DeclareUnicodeCharacter{03C9}{$\omega$}

\DeclareUnicodeCharacter{2218}{$\circ$}

\DeclareUnicodeCharacter{207B}{$^{-1}$}
\numberwithin{equation}{section}
\newtheorem{defn}[equation]{Definition}
\newtheorem{proposition}[equation]{Proposition}
\newtheorem{corollary}[equation]{Corollary}
\newtheorem{rem}[equation]{Remark}
\newtheorem{exm}[equation]{Example}
\newtheorem{lemma}[equation]{Lemma}
\newtheorem{theorem}[equation]{Theorem}
\newtheorem{notat}[equation]{Notation}
\newtheorem{newpar}[equation]{}
\newtheorem{xdefn}{Definition.}
\newtheorem{xproposition}{Proposition.}
\newtheorem{xcorollary}{Corollary.}
\newtheorem{xrem}{Remark.}
\newtheorem{xexm}{Example.}
\newtheorem{xlemma}{Lemma.}
\newtheorem{xtheorem}{Theorem.}
\newtheorem{xnotat}{Notation.}
\newtheorem{xnewpar}{\it}
\newtheorem{xproof}{{\it Proof. }}
\newtheorem{xproofof}{{\it Proof}}

\newenvironment{definition}{\begin{defn}\em}{\end{defn}}
\newenvironment{remark}{\begin{rem}\em}{\end{rem}}
\newenvironment{example}{\begin{exm}\em}{\end{exm}}
\newenvironment{notation}{\begin{notat}\em}{\end{notat}}
\newenvironment{proof}{\begin{xproof}\em}{\end{xproof}}

\newenvironment{newparagraph*}[1]{\begin{xnewpar}\hspace*{-1.5mm}{#1}. \rm}{\end{xnewpar}}
\newenvironment{definition*}{\begin{xdefn}\em}{\end{xdefn}}
\newenvironment{remark*}{\begin{xrem}\em}{\end{xrem}}
\newenvironment{example*}{\begin{xexm}\em}{\end{xexm}}
\newenvironment{notation*}{\begin{xnotat}\em}{\end{xnotat}}
\newenvironment{proposition*}{\begin{xproposition}}{\end{xproposition}}
\newenvironment{corollary*}{\begin{xcorollary}}{\end{xcorollary}}
\newenvironment{lemma*}{\begin{xlemma}}{\end{xlemma}}
\newenvironment{theorem*}{\begin{xtheorem}}{\end{xtheorem}}
\pgfplotsset{compat=1.18}

\titleformat*{\section}{\large\bfseries}

\begin{document}
\author{Clarence Protin\footnote{Centro de Filosofia da Universidade de Lisboa}}
	
	\title{Natural Term Logic}
	\maketitle 
\begin{abstract}
In this paper we develop a formal system called Natural Term Logic (NTL). NTL aims to represent key aspects of the logical and grammatical mechanisms of natural language as well as grammatical transformations which preserve core logical meaning. NTL can be seen as a refinement  of the ideas of Quine's paper `Variables Explained Away' and the technical concepts introduced by Bealer and Zalta. NTL is more fine-grained than Bealer's first-order intensional logic (BL): there is a many-to-one correspondence  $\nu$  between NTL terms and  closed BL terms as well as a canonical map  $\beta$ which assigns to each closed BL term a corresponding NTL term. The map $\nu$ can be seen as assigning a core logical content of the NTL term.
    We define a series of reductions on NTL terms which intuitivelyy speaking capture meaning-preserving syntactic transformations ( transformations which preserved the basic logical meaning of a term) and our main result is that each NTL term $T$ reduces to a unique normal term $N$. The reductions fall into the structural, predicative and pushing-in categories.  Predicative reductions decompose NTL terms so that predication is only applied to a primitive term (such terms are called prenormal). A key ingredient in the proof is the fact that $\beta \nu N = N$ when $N$ is normal. This suggests that within NTL the normal form of a term expresses the core logical content of the term.
    
\end{abstract}
 \emph{(...) for it is because of the interweaving of the forms with one another that we come to have the 
\emph{logos}}.\\

- Plato, Sophist 259E

\section{Introduction}

In this paper we develop a formal system called Natural Term Logic (NTL). NTL aims to represent key aspects of the logical and grammatical mechanisms of natural language as well as grammatical transformations which preserve core logical meaning. For quantification and multiple generality NTL makes no use of variables just as in natural language. For example: `All Alice's friends know a friend of Bob's'.  Key examples of logical meaning invariance under grammatical transformations are the following:

  \begin{enumerate}[($i$)]
    \item `Alice likes Bob' = `Bob is liked by Alice'.
  
    \item `Alice likes Alice' =  `Alice likes herself'.
    
    \item  `Alice knows that Bob sleeps' = `Bob is known by Alice to sleep'.  
   
  \end{enumerate}

Item $(i)$ is an instance of the `passive voice', item $(ii)$ an instance of the use of a ` reflexive pronoun', and item $(iii)$ an instance of `passive raising'. The latter involves what is commonly called `intensional logic' (since expressions representing propositions are subjects or objects of a predicate) but the grammatical transformation would seem to involve a complex kind of property-formation.

NTL expressions are built up inductively from primitive \emph{terms} and \emph{constructors} which take one or more terms and yield a new term.  Each term has a \emph{valency} which is a non-negative integer that represents its ``degree of saturation", $0$ for  propositions and other kinds of objects ($0$-ary relations), $1$ for properties (unary relations), $2$ for binary relations and so forth.  

Our approach is somewhat different from traditional approaches to the formalization of natural language because it is not based on variables and their bindings by operators (see e.g.\ \cite{Church,Montague}). In fact, in NTL there are no variables nor, in particular, variable binding. NTL differs from Schönfinkel and Curry's Combinatory Logic \cite{Curry} in that constructors are not terms in their own right. The precursors of NTL are the constructors used in Quine's paper `Variables Explained Away' \cite{Quine(60)} for first-order logic and the syntactic and semantic functions used by Bealer in \cite{QC} (and also used to a certain extent by Zalta in his book `Axiomatic Metaphysics' \cite{Zalta}). 

There are three families of constructors $\Upsilon, \Pi$ and $\Lambda$ which are indexed with finite combinatorial information (see Section~\ref{S:NTL} for the details). The $\Upsilon$ family can be thought of as permuting the arguments of a term, as in the construction of the passive voice in $(i)$ of `Bob is liked by Alice', where a constructor $\Upsilon$ is applied to the term representing `likes' to yield `is liked by'. The $\Upsilon$ constructor represents also a generalized diagonalization of the arguments of a term, and in particular \emph{reflexivization} of a term: applying a certain $\Upsilon$ to `likes' yields the 1-valent term `likes oneself'.  The $\Pi$ family in its simplest form expresses ordinary predication, the application of a verb to its subject and object(s): for instance the Subject-Verb-Object construction `Alice likes Bob' corresponds to the $\Pi^{(0,0) }$ application of `likes' to the arguments `Alice' and `Bob'.  But $\Pi$-applications also have more complex forms of term-formation involving `embedded predication'.  Consider the simple predication `Alice knows that Bob sleeps' which corresponds to applying the $\Pi$ constructor to  `knows',  to 'Alice' and the 0-valent term (proposition) `Bob sleeps',  itself the simple predication of `sleeps' to `Bob'. Thus `Bob' might be called an `embedded subject' which we can abstract to transform the proposition `Alice knows that Bob sleeps' into a 1-valent term (property): `the property of being known by Alice to sleep' or more simply `known by Alice to sleep'. Our constructor $\Pi$ can express the property `known by Alice to sleep' from the terms `known',  `Alice' and `sleeps'.  Herein `sleeps' does not enter directly as an argument (as in `Alice knows the property `sleeps'') but as something that contributes through its predicability to an argument to the forming a more complex property via the predication of `knows'. Having this kind of predication also allows us to express the analogue of variable substitutions (useful for the NTL versions of quantifier rules).

Finally, the family $\Lambda$ includes logical connectives as well as a way of  expressing quantification without variables or variable binding. Although we focus in the present work only on the classical connectives and universal quanifier, these logical constructors are completely general and thus can accommodate both classical and non-classical logics.

The examples $(i)-(iii)$ above can be expressed in NTL as follows:

\begin{enumerate}[$($I$)$]
\item  $\Pi^{(0,0)}$[likes][Alice][Bob] = $\Pi^{(0,0)} \Upsilon^{\{2\}\{1\}}$[likes][Bob][Alice]
\item $\Pi^{(0,0)}$[likes][Alice][Alice] = $\Pi^{(0,0)} \Upsilon^{\{1,2\}}$[likes][Alice]
\item  $\Pi^{(0,0)}$[knows][Alice]$\Pi^{(0)}$[sleeps][Bob] = $\Pi^{(0)}\Pi^{(0,1)}$[know][Alice][sleeps][Bob]
\end{enumerate}

and `all Alice's friends know Bob' can be expressed as\footnote{for simplicity we assume we have the constructor $\Lambda_\rightarrow$ which in the approach of this paper can easily be defined in terms of the logical constructors for $\neg$ and $\wedge$.}
\[
\Lambda_\forall^1\Upsilon^{\{1,2\}} \Lambda_{\rightarrow}\text{[being Alice's friend]}\Pi^{(1,0)}\text{[know]I[Bob]}.
\]

The significance of the constant $I$ will be explained further ahead. 

Each constructor of the three families is indexed with finite combinatorial information. This combinatorial information crucially depends on a theory we call \emph{weaver theory} which can be described roughly as a certain presentation and application of the theory of surjective morphisms of linearly ordered sets. Weaver theory can be interpreted geometrically as a generalization of braids in which strands can be merged together.

 In the examples $(i)-(iii)$ above not only do we intuitively assign the same basic logical meaning to the two sides but also the left-hand side is intuitively taken to be more simple or basic in some sense (the analogue of the `normal form' used in type theory): one could say that the left-hand sides render explicit the logical content implicit in the right-hand sides. Formalizing this intuition and setting up the corresponding reduction theory for NTL terms is the main goal of this paper.

The structure of this paper is as follows. In section 2 we present the preliminary technical notions and results of weaver theory required for the rest of the paper. In section 3 we present a variant of Bealer's logic and define the basic syntactic operations on terms. In section 4, which is the heart of the paper, we define Natural Term Logic (NTL) and the transformations $\beta$ and $\nu$. In section 4.1 we define the reductions for NTL and show their soundness. In section 4.2 we prove the main result of this paper concerning the normalization of NTL terms. In the conclusion we overview the results achieved and their possible applications and point out directions for future work. In the appendix we sketch a graph calculus for NTL and present some observations relating to the connection between NTL and operads.

\section{Weaver Theory}

\subsection{Weavers}

 In what follows if $A$ is a sequence we denote its $i$th element (for $i$ a non-negative integer) by $A_i$ and its \emph{length} by $|A|$. We denote the empty sequence by $\epsilon$ and consider $|\epsilon| = 0$. Given a sequence $A$ we denote by $set(A)$ the set of elements occurring in $A$. Given a non-repeating sequence $A$ we can also consider $A$ to be $set(A)$ with the linear order induced by the sequential order of $A$. 

\begin{definition} A weaver $W$ consists in a pair $(X,Y_R)$ where $X$ is a finite linearly ordered set and $Y_R$ is a linear ordering of the set of equivalence classes of an equivalence relation $R$ on $X$. \end{definition}

If the equivalence relation is clear from the context we may write $Y$ instead of $Y_R$. The equivalence relation of a weaver can be defined as $R = set(Y)$. It is straightforward to see that a weaver is equivalent to giving a surjective map of $X$ to some linearly ordered set $Y$. Also, as illustrated in the example below, one should think of $R$ as embodying the information relating to joining together of threads (generalized diagonalization) while the order on $Y$ captures the permutation of threads. 

\begin{example}\label{Ex:weaver_simple} Some simple examples of weavers are:
\begin{enumerate}[$(i)$]
    \item The \emph{identity weaver} $((x_1,...,x_n), (\{x_1\},...,\{x_n\}))$.
    \item  The \emph{empty weaver} $(\epsilon, \epsilon)$.
    \item  Let $W = ((1,2,3,4,5), (\{1,3 \}, \{2,4\},\{5\}))$. This weaver identifies threads 1 and 3, 2 and 4 into a single strand and leaves 5 untouched. This can be illustrated by the following diagram

\begin{center}
\resizebox{!}{4.5cm}{
\begin{tikzpicture}

\draw (0,0.5) node{1};
\draw (1,0.5) node{2};
\draw (2,0.5) node{3};
\draw (3,0.5) node{4};
\draw (4,0.5) node{5};
  \draw[thick] plot[smooth] coordinates {(0,0) (0,-3) (0,-5)}; 
 \draw[thick] plot[smooth] coordinates {(2,0) (0.4,-1.5) (0,-3.5)}; 
\draw[thick] plot[smooth] coordinates {(1,0) (1,-5)};
\draw[thick] plot[smooth] coordinates {(3,0) (1.4,-1.5) (1,-3.5)}; 
\draw[thick] plot[smooth] coordinates {(4,0) (2,-2) (2,-5)}; 

\draw (0, -5.5) node{\{1,3\}};
\draw (1, -5.5) node{\{2,4\}};
\draw (2,- 5.5) node{\{5\}}; 
\end{tikzpicture}}

\end{center}
\end{enumerate}
\end{example}

\begin{notation} Given a weaver $W = (X,Y_R)$ we use the notation $W^\circ = X$, $W^\flat = Y_R$ and $W_R = set(Y)$. 
\end{notation}

\begin{definition} We say that a weaver $W = (X,Y)$ has \emph{order} $(m,n)$ if $|X| = m$ and $|Y| = n$.
\end{definition}

For instance, the weaver from Example~\ref{Ex:weaver_simple}$(iii)$ has order $(5,3)$.

\begin{definition} An equivalence relation $R$ on a set $X$ is called \emph{discrete} if $xRy$ iff $x = y$. A weaver $W$ is called a \emph{permutator} if $W_R$ is the discrete relation on $W^\circ$.
\end{definition}

Note that the identity weaver is a particular case of a permutator.

\begin{definition} Let $W = (X,Y)$ be a weaver. Then we define the function
\[ W^\star : set(X) \rightarrow set(Y)\] 
as follows. Given $x \in X$ we set $W^{\star}(x) = r$ where $r$ is the element in $W_R$ such that $x \in r$.
\end{definition}

Let $M = (M_1,...,M_n)$ where the $M_i$ are sets of sets. Then we use the notation

\[ M^\cup = \left(\bigcup M_1,...,\bigcup M_n\right)\]

Analogously we define $ N^\cup$ if $N$ is a set (rather than an ordered sequence) of sets of sets.

\begin{definition} Given a weaver $W = (X,Y)$ with $Y = (Y_1,...,Y_n)$ and a linearly ordered set $A$ of the same cardinality as $X$ we denote by $W^A$ the weaver $(A,Y')$ (called the \emph{shift} of $W$ relative to $A$) defined as follows. We set $A_i(W^A)_RA_j$ iff $X_iW_RX_j$ and
we define $Y' = (Y'_1,...,Y'_n)$ where $Y'_k = \{ A_l: X_l \in Y_k\}$.
\end{definition}

\begin{definition} Let $V = (X,Y)$ be a weaver of order $(n,m)$ and $W = (Z,U)$ a weaver of order $(m,k)$. Let $W' = W^{V^\flat} = (V^\flat, U')$. The \emph{composition} of $V$ and $W$ is defined as the weaver $W'V = (X, (U')^\cup$).
\end{definition}

\begin{proposition} With the notation of the previous definition we have
\[ (WV)^\star(x) = \cup ( (W')^\star \circ  V^\star (x)) \]
\end{proposition}

\begin{definition} Given a weaver $W = (X,Y)$ we denote by $W^\sharp$ the sequence  obtained by substituting each $x_i$ in $X$ by the $r \in W_R$ such that $x_i \in R$.
\end{definition}

\begin{example}
Let $W = ((1,2,3), (\{1,3\}, \{2\}))$ then $W^\sharp = (\{1,3\}, \{2\}, \{1,3\})$.
\end{example}

The following notion will be used when we define reductions in section 4.

\begin{definition} Let $Z$ be a non-repeating sequence of length $m$ and $W = (X,Y)$ a weaver of order $(n,m)$. Then we define $W^\sharp Z$ as being obtained by replacing each $r$ in $W^\sharp$ with $Z_i$ if $r$ occurs in position $i$ in $Y$.
\end{definition}




\begin{lemma}\label{L:canonical} Given a sequence $X$ of length $n$ (possibly with repeated elements) and a non-repeating sequence $S$ such that $set(S) = set(X)$, there is a weaver $W$ such that $W^\sharp S = X$ .
\end{lemma}
\begin{proof} Let $|X| = n$ and consider the sequence $N = (1,2,...,n)$. We define the equivalence relation $R$ on $N$ by $iRj$ iff $X_i = X_j$. Note that we have a bijection $f:S \rightarrow R$ where $f(x) = r$ iff for a $n \in r$ we have that $X_n = x$\footnote{the inverse map could be defined by associating to each $r$ in $R$ the element $S_j$ where $j$ is the smallest number in $r$.}. Then $f$ determines an order on $R$ by setting $r<s$ iff $f^{-1}(r)$ comes before $f^{-1}(s)$ in $S$. Let $M$ be $R$ endowed with this linear order. It is then it is easy to check that the weaver $(N,M)$ satisfies the conditions of the lemma.
\end{proof}

\begin{notation}
We call the weaver constructed in the proof of Lemma~\ref{L:canonical} the \emph{canonical weaver} associated to $X$ and $S$ and denote it by $\mathcal{W}(X,S)$. Thus we have $(\mathcal{W}(X,S))^\sharp S = X$.
\end{notation}

The following notion is fundamental for the next subsection.

\begin{definition} Let $W =( X, Y)$ and $W' = (X', Y')$ be weavers and assume without loss of generality that $X$ and $X'$ have no common elements. Then we define $W + W' = (X + X', Y + Y')$ which is easily seen to be a weaver with $(W+ W')_R = W_R \cup W'_R$.
\end{definition}
The finite sum of weavers $\sum_{i=1}^n W_i$ can be defined in a similar way.

\begin{remark} In section 4 we consider weavers modulo shifting. Thus without loss of generality in this situation we can assume for every weaver $(X,Y)$ we have that $X = (1,2,...,n)$ for some $n$. It is then easy to see that the weaver is completely determined by $Y$ alone. We mentioned that giving a weaver of order $(m,n)$ modulo shifting is the same as giving a surjective map $(1,...,m) \rightarrow (1,...,n)$ for $n\leq m$. Thus there are exactly $n!S(m,n)$ weavers of order $(m,n)$ where $S(m,n)$ is the Stirling number of the second kind. There is no great practical advantage in presenting weaver theory in terms of surjective maps between finite linearly ordered sets. A standard way to specify such a map $f$ would require a triple $(A, B,C)$ in which $A$ is a finite sequene, $B$ is the sequence of images $f(x)$ for $x$ in $A$ and $C$ is image of $f$ as a sequence (or an order on the image). Then the definition of composition is slightly more cumbersome. The present version of weaver theory lends itself better to a geometric interpretation as well as software implementations. We sketch in the appendix a graph-theoretic notation for Natural Term Logic and the present version of weaver theory would seem to be a more intuitive basis for such a project.
\end{remark}

\subsection{Weavers and Order Partitions}

\begin{definition} Given a finite sequence $X$ (or equivalently, a finite linearly ordered set $X$ with order $<$), an \emph{order partition} $P$ on $X$ is an ordered sequence $(I_1,...,I_k)$ where each $I_i$ is either $\epsilon$ or a segment of $X$ of the form $[a_i,b_i] = \{x \in X: a_i \leq x \leq b_i\}$. The sequence of the $I_i$ must satisfy the following conditions. Let $(I'_1,...I'_m)$ be the sequence obtained by removing occurrences of $\epsilon$. Then
\begin{enumerate}
\item[$(i)$] $X = \bigcup_{i = 1}^k set(I'_i)$
\item[$(ii)$] $a'_1$ is the minimum element of $X$,  $b'_i < a'_{i+1}$ for $i < m$ and $b'_m$ is the maximum element of $X_<$.
\end{enumerate}
\end{definition}

\begin{definition} Let $W = (X,Y)$ be a weaver and $A \subset set(X)$. Then we define the \emph{restriction} of $W$ to $A$, denoted by $W_A$, as the weaver $W_X = (A,Y')$ where $A$ is considered with the induced order from $X$ and  for $Y = (Y_1,...,Y_n)$ we have that $Y'$ is obtained from $(Y_1 \cap A,...,Y_n \cap A)$ by deleting occurences of $\emptyset$. 
\end{definition}

\begin{definition} Given a sequence $A$,  an order partition $P = (I_1,...,I_k)$ on $A$ and a weaver $W$ we say that $W$ is \emph{within} $P$ if $W = W_{I_1} + ... + W_{I_k}$.
\end{definition}

The following lemma is an immediate consequence of the previous definition: 

\begin{lemma} Let $W$ be a weaver within a partition $P = (I_1,...,I_k)$ on $A$. Then $W$ induces an order partition $P' = (J_1,...,J_k)$ on $W^\flat$. 
\end{lemma}

\begin{definition} Given a weaver $W$ and a partition $P$ on $A$ we say that $W$ is \emph{without} $P$ if $W$ is not the identity weaver and its restriction $W_{I_i}$ to each $I_i$ of $P$ is the identity weaver.
\end{definition}

The main result (and algorithm) of this section is the following.

\begin{theorem}\label{T:in_out} Let $W= (X,Y)$ be a weaver and $P$ an order partition on $X$ then either $W$ is within $P$ or $W$ can be factored uniquely as $W = W_{out}W_{in}$ where $W_{in}$ is a weaver within $P$ and $W_{out}$ is a weaver without the partition $P'$ induced by $W_{in}$. 
\end{theorem}

\begin{proof} 

Let $P = (I_1,...,I_n)$. Then we define $W_{in} = \Sigma_{i=1}^n W_{I_n}$ which is evidently within $P$.
We then construct the weaver $W_{out} = ((W_{in})^\flat, Z)$ as follows. We take $Y$ and we replace each $Y_i$ with a set $\{r_1,..,r_k\}$ where the $r_i \in (W_{in})_R$ are those such that $r_i \subset Y_i$. By the definition of $W_{in}$ and the fact that $P$ is an order partition it follows that $Y_i = \cup_{l=1}^k r_i$. It is then evident that $W = W_{out}W_{in}$. To show that the factorization is unique observe that given a factorization satisfying the conditions of the theorem $W = W_{out}W_{in}$ we have that $W_{I_k} = W_{in\, I_k}$. for each $I_k$ in $P$.  If we have a second factoriation satisfying the conditions of the theorem  $W = W'_{out}W'_{in}$ then from this remark it follows that $W'_{in\, I_k} = W_{in\, I_k}$ for every $I_k$ in $P$ and hence by definition of $W_{in}$ and $W'_{in}$ being within $P$ it follows that $W_{in} = W'_{in}$.  It the follows that $W_{out} = W'_{out}$.

\end{proof}

The following simple observation plays a key role in the proof of lemma 4.21.

\begin{lemma} Let $S = (S_1,...,S_k)$ be a sequence of finite sets $S_i$ and $A$ be a non-repeating sequence such that $set(A) = set(\cup_{i=1}^n S_i)$.  For $i=1,...,n$ let $B_i$ be the elements of $S_i$ arranged in the order in which they appear in $A$. Let $C = B_1 +...+ B_n$ and $P$ be the order partition on $C$ defined by this decomposition. Let $W = (X,Y) = \mathcal{W}(C,A)$ be the canonical weaver associated to $C$ and $A$ and let $P$ be considered as defining a partition on $X$.  Then we have that $W$ is without the partition $P$.
\end{lemma}

\subsection{Actions of Number Sequences on Weavers}

\begin{definition} Given a non-negative integer $n$  and a symbol $a$ we define  for $n> 0$ $\sqcup_n a = a^1a^2...a^n$ (considered as a sequence of distinct elements all whose elements are tagged by $a$) and for $n = 0$ we set $\sqcup_0 = \epsilon$.
\end{definition}

\begin{definition} Let $W = (A,L)$ be a weaver  of order $(n,m)$ and $S$ a sequence of non-negative integers of length $m$. Then we define $W \odot S$ to be the weaver  $(A',L')$ obtained as follows.
Suppose $S_i = p$. Then we replace in $A$ each  $c \in L_i$ with $\sqcup_{p}c$. and in $L$ we replace $L_i = \{c_1,..,c_k\} $ with $L_i^1....L_i^p$ where $L_i^s = \{c^s_1,..,c^s_k\}$.  We do this simultaneously for $i=1,...m$.

\end{definition}

\begin{example} Let $L = ((1,2,3,4), (\{1,3\},\{2,4\})$ and $S = (2,2)$. Then \[W\odot S = ((1^1, 1^2, 2^1, 2^2, 3^1, 3^2, 4^1, 4^2), ( \{1^1, 3^1\}, \{1^2, 3^2\}, \{2^1, 4^1\} , \{2^2, 4^2\})\]
\end{example}

\begin{remark} In the applications that follow we will consider without loss of generality that $W \odot S$ is shifted to be defined on a sequence of the form $(1,2,...,n)$ (cf. Remark 2.16).
\end{remark}

\begin{remark} There is an alternative way of presenting $ W\odot S$.  In the conditions of the above definition let $Y = Y_1,...,Y_m$ be a non-repeating sequence of non-repeating sequences with $|Y_j| = S_j$. Consider $ Z = (W^\sharp (Y_1,...,Y_m)) ! $ where we use the notation $(X_1,...,X_m)! = X_1 + ... +X_m$ for a sequence of sequences. Then $W \odot S$ is (modulo shifting) the canonical weaver $\mathcal{W}(Z,Y)$.

\end{remark}

\begin{definition} A \emph{selector} is a pair $(N,M)$ of sequences of non-negative integers such that $|N| = |M|$ and $M_i \leq N_i$.
\end{definition}

Given a sequence of integers $N$ we denote by $\Sigma N$ the sum of the elements of $N$.  Given a symbol $s$ we denote by $s^n$ for $n$ a non-negative integer the sequence consisting of $s$ repeated $n$ times. If $n=0$ we set $s^0 = \epsilon$.

\begin{definition} Let $N= (n_1,...,n_k)$ be a sequence of non-negative integers and $S$ a sequence of length $\Sigma N$. Then the \emph{division} of $S$ by $N$, denoted by $S/N$, is the unique decomposition 

\[ S = A_1 + ...+ A_k\]

such that $|A_i| = n_k$. Note that $S/N$ is, strictly speaking, the sequence $(A_1,...,A_k)$.
\end{definition}
Note that $S/N$ defines an order partition on $S$.

\begin{definition} Given a selector $(N,S)$ with $|S| = k$ and sequence of non-negative integers $T$  with $|T| = \Sigma S$,  we define the \emph{associator} $S\otimes T$ to be a tuple $(A,B_1,...,B_k)$ constructed as follows.
Let 

\[T/S = \overline{\tau}_1 + ... + \overline{\tau}_k\]

Then we define $B_i = \overline{\tau}_i + 1^{N_i - S_i}$ for $i = 1,...,k$ and

\[A =  (\Sigma \overline{\tau}_1,...,\Sigma\overline{\tau_k})\]
\end{definition}

We use the notation $(S\otimes T)^1 = A$ and $(S\otimes T)_i = B_i$.

\begin{example}  Let $S$ be $((3,3),(2,1))$ and $T$ be $(2,2,1)$. Then  $(2,2,1) / (2,1) = (2,2) + (1)$. So $(S \otimes T)_1$ and $(S \otimes T)_2$ are the $ (2,2,1)$ and $(1,1,1))$ respectively (we have $3-2 = 1$ and $3-1 = 2$) and $(S \otimes T)^1$ is $(4,1)$. 
\end{example}


Selectors act on sequences of sequences as follows. 

\begin{definition}
Let $S = (N,M)$ be a selector with $|N| = n$ and $A = (B^1,...,B^n)$ a sequence of sequences $B^i$ with $|B^i| = N_i$. Let $C^i$ be the subsequence of $M_i$ consisting of the first $M_i$ elements of $B^i$ and let $D^i$ be the 'remainder' of $B^i$ (that is, $B^i = C^i + D^i$).
Then we set $S\bullet A = (C^1,...,C^n)$ and $A\bullet S = (D^1,...,D^n)$.
\end{definition}




\section{Bealer Logic}

Bealer logic (BL) is first-order logic with a term-forming operator $\lambda x_1...x_n \phi$ where $x_1,\dots,x_n$ is a (possibly empty) list of distinct variables and $\phi$ is a formula (or else a variable or a constant).

The language of BL consists of a countable collection of basic variables $v_1, v_2,...$, constants $c_1,c_2,...$, 
and atomic predicates $A_1,A_2,...$ of various arities. 

Variables are defined as follows:

\begin{itemize}
\item Basic variables are variables.
\item If $x_1,...,x_n$ are variables then the finite set $\{x_1,...,x_n\}$ is also a variable.
\end{itemize}

The variable  $\{x_1,...,x_n\}$ is to be considered a distinct symbol in its own right. This contruction is for book-keeping purposes and the variables that appear within the brackets  are never bound or participate in substitutions.

\begin{definition}
\emph{Formulas} and \emph{terms} are defined by mutual recursion as follows.
\begin{itemize}
\item[-] Variables and constants are terms.
\item[-] For a variable $x$ and constant $c$ the expressions $\lambda x x$, $\lambda  x$ and $\lambda  c$ are terms.

\item[-] If $t_1,...t_n$ are terms then $A(t_1,...,t_n)$ is a formula whenever $A$ is an $n$-ary atomic predicate.
\item[-] If $\phi$ is a formula and $x_1,...,x_n$ is a possibly empty list of distinct variables then $\lambda x_1...x_n\phi$ is a term. 
\item[-] If $\phi$ and $\psi$ are formulas and $x$ is a variable then $\phi\land\psi$, $\neg\phi$ and $\forall x\phi$ are formulas.
\end{itemize}
\end{definition}
Terms of the form $\lambda x_1,\dots x_n \phi$ are called \emph{complex} terms. In this case $n$ is called the \emph{sort} of the term. We consider that $\lambda x x$  has sort 1 and $\lambda x$ and $\lambda  c$ have sort $0$.

We consider equivalence modulo renaming  bound variables and discarding of dummy variables. Thus there is modulo renaming bound variables only one term of the form $\lambda x x$ (this will correspond to the 'unit' $I$ in NTL). We consider that $\lambda  x$ and $x$ are equivalent as well as $\lambda  c$ and $c$. This is convenient for the
definition of the syntactic operators further ahead. We now consider that all variables and constants are in their equivalent $\lambda$-forms.



Let $S$ be a selector and $W$ a weaver. We define the syntacitc operators $\pi^S$, $\upsilon^S$, $\lambda_\neg$, $\lambda_\wedge$ and $\lambda_\forall$ which act on BL terms to produce new terms.

\begin{definition}Let $p = \lambda x_1...x_n \phi(x_1,...,x_n)$ be a complex term and $t_1,...,t_n$ be terms written  $t_i = \lambda x^i_1...x^i_{m_i} \phi_i$ (which we write $\lambda \overline{x}^i \phi_i$ and rename the bound variables to be distinct for different values of $i$) which is of sort $m_i$. Consider a selector $S$ of the form $S = ((m_1,...,m_n), (a_1,...,a_n))$. 
 Then

\[\pi^S p t_1...t_n = \lambda  (S \bullet (\overline{x}^1,...,\overline{x}^n)) !  \phi[t'_1...t'_n/x_1...x_n] \]

where $t'_i = \lambda ((\overline{x}^1,...,\overline{x}^n)\bullet S)_i   \phi_i$.

\end{definition}
We consider that $\pi^S$ is not defined for $\lambda xx$ occurring as a head as in $\pi^A (\lambda xx )T$.

\begin{definition} let $W$ be a weaver of order $(n,m)$. Then

\[ \upsilon^W \lambda x_1...x_n \phi = \lambda (W^{(x_1,...,x_n)})^{\flat}\phi[(W^{(x_1,...,x_n)})^{\sharp}  / x_1...x_n   ]  \]

\end{definition}

\begin{definition}
Let $p = \lambda x_1...x_n \phi(x_1,...,x_n)$ and $q = \lambda y_1...y_m \psi(x_1,...,x_m)$ be complex terms (where if necessary we rename the bound variable so that the $x_i$ and $y_j$ are distinct). Then we define 

\begin{align*}
\lambda^n_\neg p &= \lambda x_1...x_n \neg\phi(x_1,...,x_n)\\
\lambda^{n,m}_\wedge p q &= \lambda x_1...x_ny_1...y_m \phi(x_1,...,x_n) \wedge \psi(y_1,...,y_m) \\
 \lambda^{n}_\forall p &= \lambda x_2...x_n \forall x_1\phi(x_1,...,x_n).
\end{align*}

\end{definition}

\begin{remark} The operator $\pi$ includes, as a particular case, standard substitution of a variable for a term. Also note that we have a trivial instance of the application  $\pi^{(1,1,...,1)}t(\lambda x_1 x_1)...(\lambda x_n x_n)$ which yields $t$ itself.
\end{remark}

\begin{definition}
A \emph{primitive} BL term is a term of the form $\lambda x_1...x_n A(x_1,...,x_n)$ where $A$ is an atomic predicate. We also use the notation $[A]$. 
\end{definition}

It will follow from the results of the next section that any closed complex term $t$ in BL can be obtained by successive applications of $\pi, \sigma, \lambda_\ast$ for $(\ast \in \{\neg, \wedge, \forall \})$ to primitive terms corresponding to the atomic predicates occurring in $t$.

\begin{example}
The term  $\lambda w A(\lambda v B(w,v), \lambda u C(u))$  is obtained from  $\pi^{(1,0)} [A][B][C]$.
\end{example}



\begin{definition}
A term of the form $\lambda x_1...x_m A(t_1,...,t_n)$ is called an \emph{atomic} term.
\end{definition}

\section{Natural Term Logic}\label{S:NTL}

The language of Natural Term Logic (abbreviated NTL from now on) is defined as follows. We are given a collection of \emph{primitive} terms $A_1$, $A_2$,...,$A_n$,... including the special term $I$ and a finite collection of constructors $\Gamma_1,\Gamma_2,...,\Gamma_n$. Each primitive term $A$ has a \emph{valency} (saturation degree) which is a non-negative integer $s$. 
 Each constructor $\Gamma$ has a signature $(s_1,...,s_n)\rightarrow s_{n+1}$ where the $s_i$ are non-negative integers.

\begin{definition}
\emph{Terms} in NTL are defined as follows
\begin{itemize}
\item[-] A primitive term is a term.
\item[-] If $\Gamma$ is a constructor with signature $(s_1,...,s_n)\rightarrow s_{n+1}$ and $T_i$ are terms with valencies $s_i$ for $i=1,...,n$ then $\Gamma T_1...T_n$ is a term of valency $s_{n+1}$.
\item[-] Nothing else is a term.
\end{itemize}
\end{definition}

We sometimes indicate that a term $T$ has valency $s$  by writing $T^{(s)}$. The special term $I$ has valency $1$.

In this paper we assume the following relationship between the language of BL and that of NTL. For every $n$-ary atomic predicate symbol $A$ of BL there is a unique corresponding primitive term $A^{(n)}$ of NTL and vice-versa. We can also assume that NTL has a set of primitive terms of valency 0 which correspond to constants in BL. We will also work only with the following set of constructors:

\[ \Pi^S, \Upsilon^W,  \Lambda^{(n,m)}_\wedge, \Lambda^n_\neg, \Lambda^n_\forall\]
where $S$ is a selector and $W$ is a weaver defined on a (possibly empty) sequence of the form $(1,2,...,n)$. We have that $\Pi^S$ has signature $(n, a_1,...,a_n) \rightarrow \Sigma b$ for $S = (a,b)$, $\Upsilon^W$ has signature $(n) \rightarrow m$ for $W$ of order $(n,m)$, $\Lambda^{(n,m)}_\wedge$ has signature $(n,m) \rightarrow n + m$ for $n,m\geq 0$,  $\Lambda^n_\neg$ has signature $(n) \rightarrow n$ and  $\Lambda^n_\forall$ signature $(n) \rightarrow n - 1$ for $n \geq 1$.
For a selector $S =(a,b)$ we will write $\Pi^SA B_1...B_n$  as $\Pi^bA B_1...B_n$ because $a$ is determined by the valencies of the $B_i$  (cf. the NTL terms discussed in the Introduction). And as mentioned before if we assume that all weavers are defined on $(1,2,...,n)$ then for $W = (X,Y)$ we need only indicate $Y$ in the superscript $\Upsilon^Y$. And importantly, we do not allow the special primitive term $I$ to occur as the head of $\P$i-application or as an argument  of the constructors $\Lambda_\neg, \Lambda_\wedge$ or $\Lambda_\forall$. Thus for instance  $\Pi^A I T$ is not a valid NTL term.

We now define a canonical process $\nu$ to convert an NTL term $T$ into a closed BL term $\nu T$. 

\begin{definition}
We define $\nu$ by induction on the structure of $T$.
\begin{itemize}
    \item[-] $\nu A^{(n)} := \lambda x_1...x_n A(x_1,...,x_n)$
    \item[-] $\nu I = \lambda x  x$
    \item[-] $\nu \Gamma T_1...T_n := \gamma (\nu T_1)...(\nu T_n),$\\ where, on the left-hand side  $\gamma$ corresponds to associating the constructor with the corresponding BL operator:
     $\Pi^S$ becomes $\pi^S$, $\Upsilon^W$ becomes $\upsilon^W$, $\Lambda^{(n,m)}_\wedge$ becomes $\lambda^{(n,m)}_\wedge$ and so forth.
    \end{itemize} 
\end{definition}

The following is clear from the definition of $\nu$.

\begin{proposition} To every term $T$ in $\rm NTL$ we can associate canonically a closed term $\nu T$ of $\mathrm{BL}$ in such a way that for a primitive term $A^{(n)}$ we have that $\nu A^{(n)} := \lambda x_1...x_n A(x_1,...,x_n)$.
\end{proposition}


\begin{definition}\label{D:Bealerdecomposition}
Given a closed BL term $t$, its \emph{Bealer decomposition} $\beta t$  is a process which associates $t$ to an NTL term $\beta t$ defined according to the following clauses: 
\begin{itemize}
\item[-] If $t$ is a primitive term $[A]$  then $\beta t =  A$ and if $t$ is a constant $c$ then $\beta t = C$ (the valency 0 primitive term in NTL corresponding to $c$).

\item[-] If $t$ is $\lambda xx$ then $\beta t = I$.

\item[-] Let $t$ be a non-primitive atomic term $\lambda x_1...x_m A(t_1,...,t_n)$ where the $t_i$ have sort $n_i$.  Let $V_i$ be the sequence of free variables in $t_i$   ordered according to the order they appear in $x_1,...,x_m$.  Consider the sequence $V = V_1 + ...+V_n$ and let $W = \mathcal{W}(V,x_1,...,x_n)$. 


 Then we define:

\[ \beta t =  \Upsilon^ W \Pi^{(|V_1|,...,|V_n|)}A (\beta \lambda V_1 t_1)...(\beta \lambda V_n t_n),  \]

By $\lambda V_k t_k$ is understood that we place the variables $U_k$ right after $\lambda$ in $t_k$. If $V_k$ is empty then this expression is just $t_k$.


If $W$ is the trivial weaver then $\Upsilon^W$ is omitted and if the $t_i$ are all variables then we set simply

\[ \beta t =  \Upsilon^ W  A   \]

\item[-] If $t$ is $\lambda x_1...x_n \phi \wedge \psi$. Let $V_1$ be the free variables in $\phi$ ordered as they occur in $x_1...._n$ and $V_2$ the free variables in $\psi$ ordered as they occur in $x_1...x_n$. Let $V = V_1 + V_2$ and $W = \mathcal{W}(V, x_1...x_n)$. 


Then we set

\[\beta t=  \Upsilon^W \Lambda^{(|V_1|, |V_2|)}_\land (\beta \lambda V_1\phi )( \beta \lambda V_2 \psi)  \]

where if $W$ is the trivial weaver we omit $\Upsilon^W$.

\item[-] If $t$ is $\lambda x_1...x_n \neg \phi$, then $\beta t = \Lambda^n_\neg (\beta \lambda x_1...x_n  \phi).$

\item[-] If $t$ is of the form $\lambda x_1...x_n \forall y \,\phi(y)$, then $\beta t =  \Lambda^{n+1}_\forall (\beta \lambda yx_1...x_n\phi).$
\end{itemize}
\end{definition}

Note that since $t$ is a closed term, in the above definition there is never a case in which we need to define $\beta \lambda x$

\begin{theorem} $\beta t$ is well-defined NTL term and we have $\nu\beta t = t$.
\end{theorem}
\begin{proof}
By a straightforward induction on the structure of $t$.
\end{proof}

\begin{corollary} A closed complex term $t$ of $\mathrm{BL}$ is always the result of the application of some composition of the syntactic operators  $\pi, \sigma, \lambda_\ast$ for $(\ast \in \{\neg, \wedge, \forall \}$ to $[A_1],...,[A_k]$ and possibly $\lambda xx$  and some constants where $A_1,...,A_k$ are the atomic predicates occurring in $t$.
\end{corollary}

Note that in general the decomposition of a BL term into compositions of the syntactic operators is not unique. Consider $t = \lambda yx A(\lambda B(x,y))$. Then this term can be decomposed either in the form $\Upsilon^W \Pi^S [A][B]$ or $\Pi^T [A] \Upsilon^V [B]$.

\begin{remark}  Bealer in \cite{QC} uses a syntactic decomposition of closed terms based on a different, more complex system, than our $\beta$ in $NTL$. Instead of purely logical constructors Bealer uses models and semantic operators on such models which can be described as a semantic counterpart of a far more complex and convoluted system of constructors. In Bealer's system it is non-trivial to show that the interpretation $\beta t(x)[s/x]$ corresponds to $\pi^{(0)} (\beta \lambda x  t) (\beta s) $, the proof of which, if it where given, would require the development of the semantic version of something similar to the NTL reductions which we  expound in the next subsection.

\end{remark}

\begin{remark} In NTL there are no primitive terms specifically for individuals unlike in \cite{QC}. Our idea is to subsume both individuals and propositions into the general notion of a primitive term of valency 0 - a term which can never be the head of a $\Pi$-application. The distinction could be determined by the property of being a truth-bearer.
\end{remark}

\subsection{Reductions}

We will define a series of reductions for NTL. The central idea is to reduce the complexity of the head of a $\Pi$ application and to organize $\Upsilon$ applications in a canonical way. Reductions will allow us to define the concepts of normalization and normal form for NTL terms.  It will follow from the main result of this paper that $T_1$ and $T_2$ reduce to the same normal term $N$ iff $\nu T_1 = \nu T_2$.

Reductions will be divided into

\begin{itemize}
\item[-]  Structural reduction: the merging of sucessive $\Upsilon$-applications and dropping of trivial $\Upsilon$-applications.
\item[-] Predicative reductions, involving the removing of complex terms from the head of $\Pi$-applications and the simplification of trivial applications to series of $I$s.
\item[-] Pushing-in reductions.
\end{itemize}

We introduce the reductions progressively while also  checking their `soundness', that is, if the associated BL term remains the same. A rule $A \rightsquigarrow B$ is called \emph{sound} if $\nu A = \nu B$. The soundness of the reduction rule plays a fundamental role in the proof of our main result about normalization in the next section.

The structural reductions consist in two rules:

\begin{equation}
\tag{$C_\Upsilon$} \Upsilon^{W_1} \Upsilon^{W_2} X \rightsquigarrow \Upsilon^{W_1W_2} X
\end{equation}

and

\begin{equation}
\tag{$Id_\Upsilon$} \Upsilon^{\{1\}...\{n\}} X^{(n)} \rightsquigarrow  X
\end{equation}

Predicative reductions have as goal to progressively simplify the head $T$ of a $\Pi$ application to the case in which the head is a primitive term. Thus we have reductions for each of the types of complex term $T$.

Let us start with the case in which $T$ begins with $\Pi$. We postulate reductions of the form:

\begin{equation}
\tag{$R_{\Pi}$}
\Pi^A(\Pi^B TT_1...T_n)S_1...S_m \rightsquigarrow \Pi^{(C\otimes A)^1} T U_1...U_n
\end{equation}

where $U_i = \Pi^{(C\otimes A)_i}T_i\overline{S'}^n$ if $T_i$ is not $I$ and $U_i = S_i$ otherwise, where $C = (V, B)$ is the selector determined by $B$ and the valencies of $T_1,...,T_n$, $S_1...S_m / B = \overline{S}^1 +...+\overline{S}^n$  and $\overline{S'}^i = \overline{S}^i I^{(V_i - B_i)} $.

Assume that $T$ is a $\Upsilon^W$ application with $W$ a weaver. Then we postulate:

\[
\tag{$R_\Upsilon$}
\Pi^{A}(\Upsilon^W T)S_1...S_n \rightsquigarrow \Upsilon^{W\odot A} \Pi^{ W^\sharp A} T W^\sharp (S_1...S_n)\]

When the head $T$ is an $\Lambda^n_\neg$ application we postulate:

\[
\tag{$R_\neg$}
\Pi^{A}(\Lambda^n_\neg T)S_1...S_n ) \rightsquigarrow \Lambda^{\Sigma A}_\neg  (\Pi^{A} T S_1...S_n).\]

The case in which $T$ is an $\Lambda^{n,m}_\wedge$ application is as follows. We postulate

\[
\tag{$R_\land$}
\Pi^{A}(\Lambda^{n,m}_\land T^{(n)}S^{(m)})T_1...T_{n+m} \rightsquigarrow \Lambda^{\Sigma A_1, \Sigma A_2}_\land (\Pi^{A_1}TT_1...T_n)(\Pi^{A_2}ST_{n+1}...T_m)  \]

Here $A /(n,m) = A_1 + A_2$.


For $\Lambda^n_\forall$ applications we postulate:

\[
\tag{$R_\forall$}
\Pi^{A}(\Lambda^{n+1}_\forall T^{(n+1)})T_1...T_n \rightsquigarrow \Lambda^{1 + \Sigma A}_\forall (\Pi^{(1) + A} TIT_1...T_n)    \]

Then finally we have the predicative reduction: 

\[
\tag{$R_{I}$}
\Pi^{((1,1,...,1)} T II....I \rightsquigarrow T\]


Suppose we have a term of the form $\Upsilon^W\Pi^A TS_1...S_n$. By Theorem\ref{T:in_out}, $W$ can be decomposed as $W_{out}W_{in}$ for $W_{in}$ within the order partition $(1,2,...,\Sigma A) / A$. Consider the decomposition induced by the partition: $W_{in} = \Sigma_i W_i$ with $W_i$ of order $(o_i, p_i)$. Let $W$ be a weaver of ordern $(n,m)$. We use the notation $W+ 1^k$ to indiciate the sum of $W$ and the identity weaver $(n+1, n+1,...,n+k),\{n+1\}...\{n+k\})$.  Then we postulate the reduction

\[\tag{$P_\Pi$} \Upsilon^W\Pi^{A} TS_1...S_n \rightsquigarrow \Upsilon^{W_{out}} \Pi^{A'}T (\Upsilon^{W_1 +1^{V_1-A_1}}S_1)...(\Upsilon^{W_n + 1^{V_n-A_n}}S_n)             \]

where $A' =    (p_1,...,p_n)$ and $V_i$ is the valency of $S_i$.


Let $W$ be a weaver of order $(n,m)$ decomposed as $W_{out}W_{in}$ for the partition  $(1,2,...,n+m) / (n,m)$ with corresponding decomposition $W_{in} = W_1 + W_2$ of orders $(n,o)$ and $(m,p)$. Then we postulate:

\[ 
\tag{$P_\land$}
\Upsilon^W\Lambda^{n,m}_\wedge TS \rightsquigarrow \Upsilon^{W_{out}} \Lambda^{o,p}_\wedge (\Upsilon^{W_1}S)(\Upsilon^{W_2}T).   \]


The case of $\Upsilon^W \Lambda^n_\neg T$ for $W$ of order $(n,m)$ is:

\[ \tag{$P_\neg$} \Upsilon^W \Lambda^n_\neg T \rightsquigarrow \Lambda^m_\neg \Upsilon^W T . \]


Finally for $\Upsilon^W \Lambda^n_\forall T$ we postulate:

\[\tag{$P_\forall$} \Upsilon^W \Lambda^n_\forall T \rightsquigarrow \Lambda^{m+1}_\forall \Upsilon^{1+W} T  \]

where $W$ has order $(n-1,m)$. 

We note that the pushing-in reductions are in some sense redundant if the $\Upsilon^W$ application is a head of a $\Pi^A$ application. For the $\Upsilon^W$ can be eliminated by a predicative reduction (note that it will follow from our main result that all reduction paths lead to the same normal term). Pushing-in reductions are important  for instance for having a unique canonical version of an outermost $\Upsilon^W$ application of a term as well as a unique form for non-head arguments of $\Pi^A$- applications.  Pushing-in reductions are essential to obtain our result concerning reduction to a unique normal term.

\begin{lemma}\label{L:C_sound} The reduction $C_\Upsilon$ and $Id_{\Upsilon}$ are sound. \end{lemma}

\begin{proof} This soundness of the first reduction amounts to showing that

\[ \upsilon^W \upsilon^V \lambda X \phi(X) = \upsilon^{WV}\lambda X  \phi(X)\]

for $V$ of order $(n,m)$, $W$ of order $(m,k)$ and $|X| = n$. We can assume $V$ and $W$ have been shifted so that $V^\circ = X$ and $W^\circ = V^\flat$. The result then follows from a direct computation of both sides of the equality followed by renaming each bound variable $\{X_1,...,X_m\}$ of the result of the left side by $\cup \{X_1,...,X_m\}$. For example, a bound variable $\{\{x,y\}, \{z\}\}$ is renamed to $\{x,y,z\}$. The soundness of $Id_{\Upsilon}$ is trivial.




\end{proof}


\begin{remark}
In what follows we will use the general notation $[X_i]_{i=1}^n$ as follows according to the context. $\lambda [X_i]_{i=1}^n \phi$ stands for $\lambda X_1...X_n \phi$ when $X$ is a list of variables (including the empty case and the case of a single variable).  If $[A_i/B_i]$ are (simultaneous) substitutions then $[A_i/ B_i]_{i=1}^n$ denotes $[A_1....A_n /B_1....B_n]$. We also use the notation $\Pi^A T[S_i]_{i=1}^n$ for $\Pi^ATS_1....S_n$.
Given a sequence $X$ of length $n$ and $m \leq n$ a non-negative integer we consider the spliting $X = X^{\leftarrow} + X^{\rightarrow}$ where $|X^{\leftarrow}| = m$.
\end{remark}

\begin{lemma}  The reduction $R_\Pi$ is sound. \end{lemma}

\begin{proof} Let $t = \nu T = \lambda X \phi(X)$, $t_i = \nu T_i = \lambda Y^i\psi_i(Y^i)$ for $i=1,...n$ and $s_j = \nu(S_j) = \lambda Z^j\rho_j(Z^j)$ for $j = 1,...m$ where we have $m = \Sigma_{i = 1}^n |Y^i| = \Sigma B$ and $|B| = n$.  We consider the splitting $Y^i = Y^{i\leftarrow} + Y^{i\rightarrow}$ determined by $B_i$ and $Z^j = Z^{j\leftarrow} + Z^{j\rightarrow}$ determined by $A_j$.

Then the computation of $\nu$ of the left side yields

\[ \pi^A\pi^B t [t_i]_{i=1}^n [s_j]_{j=1}^m  = \pi^A \lambda [Y^{i\leftarrow}]_{i=1}^n \phi [    [ \lambda Y^{i\rightarrow}\psi(Y^{i\leftarrow},Y^{i \rightarrow}) ]_{i=1}^n        / X]\]

Let

\[ [Z^{j\rightarrow}]_{j=1}^m / B= [Z'^{i^\rightarrow}]_{i=1}^n \]

\[ [Z^{j\leftarrow}]_{j=1}^m  /B = [Z'^{i^\leftarrow}]_{i=1}^n \]

and

\[[ s_j]_{j=1}^m / B  = \overline{s}^1 + ... + \overline{s}^n \]

\[[ \rho_j]_{j=1}^m / B  = \overline{\rho}^1 + ... + \overline{\rho}^n \]

Then computing the last expression above yields

\[ =   \lambda [Z'^{i\leftarrow}]_{i=1}^n \phi [    [ \lambda Y^{i\rightarrow}\psi_i(Y^{i\leftarrow},Y^{i \rightarrow})]_{i=1}^n        / X][  [\lambda Z'^{i\rightarrow}_k\overline{\rho}^i_k(Z'^{i\leftarrow},Z'^{i \rightarrow})]_{k=1}^{B_i} /Y^{i\leftarrow} ] _{i=1}^n\]

which is evidently equal to

\[ =   \lambda [Z'^{i\leftarrow}]_{i=1}^n \phi [    [ \lambda Y^{i\rightarrow}\psi_i(Y^{i\leftarrow},Y^{i \rightarrow}) [ [ \lambda Z'^{i\rightarrow}_k\overline{\rho}^i_k(Z'^{i\leftarrow},Z'^{i \rightarrow})]_{k=1}^{B_i}  /Y^{i\leftarrow} ] ]_{i=1}^n       / X] \tag{1}\]

Now consider the terms

\[ u_i =  \lambda Z'^{i\leftarrow}Y^{i\rightarrow}\psi(Y^{i\leftarrow},Y^{i \rightarrow}) [ [ \lambda Z'^{i\rightarrow}_k\overline{\rho}^i_k(Z'^{i\leftarrow},Z'^{i \rightarrow}) ]_{k=1}^{B_i}/Y^{i\leftarrow} ]\]

Then it is evident that $u_i = \nu (\Pi^{(C\otimes A)_i} T_i \overline{S}^i I^{(V_i -B_i)} ) = \pi^{(C\otimes A)_i} t_i \overline{s}^i [\lambda x,x]_{i = 1}^{V_i -B_i} $   and that the expression (1) is equal to

\[ = \pi^{(C\otimes A)^1} tu_1....u_n \]

and so the result follows.

\end{proof}

\begin{lemma}  The reduction $R_\Upsilon$ is sound.\end{lemma}
\begin{proof}

Let $t = \nu T= \lambda X\phi(X)$ and $s_i = \nu  S_i = \lambda Y^i\psi_i(Y)$ for $i=1,...,n$. Assume that $W$, of order $(|X|,n)$ is shifted to $X$. Then calculating $\nu$ of the left-hand side of $R_\Upsilon$ yields

\[ \pi^A (\lambda W^\flat \phi(X) [  W^\sharp /  X] ) [  \lambda Y^i\psi_i(Y) ]_{i=1}^n = \lambda [Y^{i\leftarrow}]_{i=1}^n \phi(X) [  W^\sharp /  X] [ [\lambda Y^{i\rightarrow} \psi_i(Y) ]_{i=1}^n  / W^\flat])\]

But this can be written as

\[ \lambda [Y^{i\leftarrow}]_{i=1}^n \phi(X) [  W^\sharp [\lambda Y^{i\rightarrow} \psi_i(Y) ]_{i=1}^n  /  X ] \tag{$\star$} \]

where we rename bound variables for different occurrences of the same $\lambda Y^{i\rightarrow} \psi_i(Y)$. Thus we could have different occurrences of the same term renamed $\lambda Y^{(1)i\rightarrow} \psi_i(Y^{(1)})$, $\lambda Y^{(2)i\rightarrow} \psi_i(Y^{(2)})$, and so forth.

We wish to put this last expression in the form

\[ \upsilon^{A'} \pi^{W'} t W^\sharp [s_i]_{i=1}^n \]

Consider the term

\[ u =  \lambda W^\sharp [Y^{i\leftarrow}]_{i=1}^n \phi(X) [  W^\sharp [\lambda Y^{i\rightarrow} \psi_i(Y^i) ]_{i=1}^n  /  X ]  \]

where we rename through numeric supscripts bound variables for different occurrences of $Y^{i \leftarrow}$ as we did above.

Then evidently $u = \nu (\Pi^{W^\sharp A}TW^\sharp [S_i]_{i=1}^n)$. The result follows from $\upsilon^{W\odot A} u = (\star)$ bearing in mind remark 2.28 and renaming variables $\{ Y^{(1)i\leftarrow}_k,  Y^{(2)i\leftarrow}_k, ...\}$ to  $Y^{i\leftarrow}_k$.

\begin{remark} Let us elucidate the above proof with a concrete example.   Let $t = \lambda xyz A(x,y,z)$, $s_1 = \lambda x_1x_2 B(x_1,x_2)$ and $s_2 = \lambda y_1y_2y_3C(y_1,y_2,y_3)$. Consider the weaver $W = ((x,y,z), \{z\}\{x,y\})$. Then $W^\sharp = \{x,y\}\{x,y\}\{z\}$ and ($\star$) corresponds in this case to 

\[\pi^{(1,2)}\upsilon^W ts_1s_2 = \pi^{(1,2)}\lambda \{z\}\{x,y\}) A( \{x,y\}, \{x,y\}, \{z\}) _1s_2 = \pi^{(1,2)}\lambda wv A(v,v,w) s_1s_2\]
\[ = \lambda x_1y_1y_2 A(\lambda y_3C(y_1,y_2,y_3), \lambda y_3C(y_1,y_2,y_3), \lambda x_2B(x_1,x_2)) \tag{a}\]

We have $W^\sharp s_1s_2 = s_2s_2s_1$ and $W^\sharp(1,2) = (2,2,1)$ so that expression $u$ in the previous proof is in this case

\[\pi^{W\sharp(1,2)} tW^\sharp s_1s_2 = \pi^{(2,2,1)} \lambda xyzA(x,y,z) s_2s_2s_1  \]
\[=  \lambda y_1y_2y'_1y'_2x_1  A(\lambda y_3C(y_1,y_2,y_3), \lambda y'_3(y'_1,y'_2,y'_3), \lambda x_2B(x_1,x_2)) \tag{b} \]

(recall how we rename bound variables for repeated applications of the term $s_2$).

Let us use remark 2.28 to calculate $W\odot (1,2)$. We can use $Y_1Y_2$ where $Y_1 = x$ and $Y_2 = y_1y_2$, corresponding to $S = (1,2)$. Then $(W^\sharp (Y_1,Y_2))! = Y_2Y_2Y_1 = y_1y_2y_1y_2x$. So $W\odot (1,2) = \mathcal{W}(xy_1y_2, y_1y_2y_1y_2x_1)$ which we can represent as 

\[ W\odot (1,2) = ((1,2,3,4,5),  \{5\}\{1,3\} \{2,4\}) \]

(one can check that $W^\sharp xy_1y_2 = y_1y_2y_1y_2x$). It is then immediate to check that $\upsilon^{W\odot S}$ of expression (b) is equal (modulo renaming of bound variables, $\{y_i,y'_i\}$ to $y_i$, for $i=1,2$ and $\{x\}$ to $x$) to expression (a) thus establishing that

\[ \pi^{(1,2)}\upsilon^W ts_1s_2  = \upsilon^{W\odot (1,2)} \pi^{W\sharp(1,2)} tW\sharp s_1s_2 \]

\end{remark}

 \end{proof}









\begin{lemma} The reductions  $R_\wedge$, $R_\neg$, $R_\forall$ and  $R_{I}$ are sound.
\end{lemma}

\begin{proof}  The soundness of these  reductions is straightforward. We show as an example the soundness of $R_\neg$.
Let $t = \lambda x_1...x_n T(x_1,...,x_n)$ and $s_i = \lambda \overline{a}^i\overline{b}^i B_i(\overline{a}^i, \overline{b}^i)$ for $i=1,...,n$ and $|\overline{a}^i | = S_i$ and $|\overline{b}^i| = N_i - S_i$. 

Then a direct calculation yields

\[\nu (\Pi^S(\Lambda^n_\neg t)s_1...s_n )) = \lambda \overline{a}^1....\overline{a}^n \neg A(x_1,...,x_n)[ \lambda \overline{b}^i B_i(\overline{a}^i,\overline{b}^i) / x_i]_{i = 1}^n\]

which is evidently equal to

\[\nu( \Lambda^{\Sigma S}_\neg  (\Pi^S t s_1...s_n)   )  \]

\end{proof}

\begin{lemma}\label{L:PushinigSoundness}
    The pushing-in reductions are sound.
\end{lemma}

\begin{proof} This follows from a straightforward computation.
\end{proof}

 




\begin{theorem}\label{T:Soundness}
The reductions of NTL are sound
\end{theorem}
\begin{proof}
    The result follows immediately from Lemmas~\ref{L:C_sound}--\ref{L:PushinigSoundness}.
\end{proof}

\subsection{Normalization}

We use the notation $T \rightsquigarrow^\ast S$ to indicate that $S$ is obtained from $T$ by multiple applications of reductions to subterms of $T$.

\begin{definition}
    An NTL term is called \emph{prenormal} if for all subterms of the form $\Pi T T_1...T_n$ we have that $T$ is a primitive term and there are no subterms of the form $\Pi^{((1,1,...,1), (1,1,...,1))} T I....I$ or $\Pi^{(1)} IT$.
\end{definition}

Note that applying only predicative reductions we always arrive in finitely many steps at a prenormal term. So prenormal terms do not allow further predicative reductions.


\begin{definition} A term of the form $\Upsilon^W T$ where $T$ starts with $\Pi$ or $\Lambda$ is called \emph{pushed-in} if the pushing-in reductions cannot be applied. 
\end{definition}

\begin{definition} A prenormal NTL term is called \emph{normal} if it all its subterms of the form   $\Upsilon^W T$ are pushed-in and the structural reductions cannot be applied. 
\end{definition}

It is evident that any sequence of reduction rules applied to subterms of a term $T$ must eventually terminate. Furthermore, the soundness of the reductions permits us to conclude directly that:

\begin{lemma} Given any NTL term $T$ any sequence of applications of the reduction rules eventually terminates in a normal term $N$, $T\rightsquigarrow^\ast N$, and we have $\nu T = \nu N$.
\end{lemma}

The main result of this paper is that such a normal term $N$ is unique. To show this we need the following lemma:


 


\begin{lemma} Let $N$ be a normal term. Then we have $\beta \nu N = N$.
\end{lemma}

\begin{proof} By induction on the structure of the normal term $N$.  If $N$ is primitive then the result follows immediately from the definition of $\beta$ and $\nu$.  Let now $N$ be a normal term.  If $N$ is of the form $\Upsilon^W P$ where $P$ is primitive then the result again follows from the definitions of $\beta$ and $\nu$.  If $N$ is of the form $\Pi^S TS_1...S_n$ then $T$ must be primitive and the result follows directly from the definitions of $\nu$ and $\beta$ and from the induction hypothesis. The cases of $\Lambda_\neg N'$ and $\Lambda_\forall N'$ are also straightforward. 

Assume now that $N = \Upsilon^W \Pi^SAB_1....B_n$. Then $A$ must be primitive and $W$ both non-trivial and without the partition defined by the selector $S$. Then we have

\[ \beta \nu \Upsilon^W \Pi^SAB_1....B_n = \beta \upsilon^W \pi^S (\nu A)(\nu B_1)...(\nu B_n) \]

Since $W$ is without the partition defined by $S$ we have that $\upsilon^W$ preserves the order of the first $S_i$ variables in $\nu B_i$ as they appear in $\pi^S (\nu A)(\nu B_1)...(\nu B_n)$.  By the definition of $\beta$ and the uniqueness of a weaver without a partition we then have

\[ \beta \upsilon^W \pi^S (\nu A)(\nu B_1)...(\nu B_n) =   \Upsilon^W \Pi^S (\beta\nu A) \beta t'_1....\beta t'_n  \]

Now $t'_i = \nu B_i$ because $W$ is without the partition defined by $S$ and so does not permute the first $S_i$ variables after the $\lambda$ (cf. lemma 2.23). Hence the result follows by induction.  Note that since $N$ is normal the head $A$ cannot be $I$ and at least one $B_i$ must not be $I$ (and note also that we assume that a variable $x$ is represented as $\lambda  x$ in $\nu N$).

An entirely analogous argument applies to the final case in which $N = \Upsilon^W\Lambda^{(n,m)}_\land AB$.

\end{proof}



\begin{remark} We illustrate the main step of the proof above through a concrete example. We take the case of a normal term $N = \Upsilon^W \Pi^SAB_1B_2$ in which $\nu(A) = \lambda xy A(x,y)$, $\nu (B_1) = \lambda x_1x_2 B_1(x_2,x_1)$ and $\nu(B_2) = \lambda y_1y_2y_2B_2(y_3,y_2,y_1)$, $S = ((2,3), (2,2))$ and $W = ((1,2,3,4), (\{1,3\}, \{2, 4\}))$ (which is clearly without the partition $((1,2), (3,4))$, and hence restricts to the identity - and so preserves the order - on the intervals of the partition\footnote{and thus $\upsilon^W$ when restricted to $x_1x_2$ and $y_1y_2$ preserves the corresponding order of the variables so that these sequences are restored when $\beta$ is calculated (cf. lemma 2.23).}). Then a direct calculation of $\nu(N)$ yields

\[\upsilon^W \pi^S (\nu A)(\nu B_1)(\nu B_2) = \lambda \{x_1,y_1\} \{x_2,y_2\} A(\lambda B_1(\{x_2,y_2\}, \{x_1,y_1\}), \lambda y_3 B_2(y_3, \{x_2,y_2\}, \{x_1,y_1\}))\]
Now consider the application of $\beta$ to the above term. By the definition of $\beta$ the weaver $W'$ for the resulting $\Upsilon$ constructor is obtained as $\mathcal{W}(\{x_1,y_1\}\{x_2,y_2\}\{x_1,y_1\}\{x_2,y_2\}, \{x_1,y_1\}\{x_2,y_2\})$ which is equal to $W$. The selector obtained from the definition of $\beta$ is clearly equal to $S$ and the terms that appear in the scope of $\beta$ after the head argument of $\Pi$ result from $\nu(B_1)$ and $\nu(B_2)$ by renaming bound variables. Thus, since we can easily check that $\beta\nu B_1 = B_1$, $\beta \nu B_2 = B_2$ and $\beta \nu A = A$ we obtain that
\[ \beta \nu \Upsilon^W \Pi^SAB_1B_2 = N \]
\end{remark}

We can now prove the main result of this paper

\begin{theorem} Let $T$ be NTL term. Then there is a unique normal term $N$ such that $T\rightsquigarrow^\ast N$.

\end{theorem}

\begin{proof} Assume that we had normal terms $N_1$ and $N_2$ such that $T \rightsquigarrow^\ast N_1$ and $T \rightsquigarrow N_2$. Then $\nu T = \nu N_1$ and $\nu T = \nu N_2$ and thus $\nu N_1 = \nu N_2$. Since $N_1$ and $N_2$ are normal we may apply the above lemma to obtain $\beta \nu N_1 = N_1$ and $\beta \nu N_2 =  N_2$. But since $\nu N_1 = \nu N_2$  the result follows.

\end{proof}

Since it can be checked by a simple inductive argument that $\beta t$ is always a normal NTL term, the above result also provides us with an alternative method to calculate the normal form of a NTL term $T$.  Instead of applying reductions directly to $T$, we simply compute $N' = \beta \nu T$.  Then since $\nu N' = \nu \beta \nu T = \nu T$ if $N$ is the normal form of $T$ then $N = \beta \nu N = \beta \nu T = N'$.

\begin{example} The author has developed a tool that allows you to interactively reduce NTL terms until reaching a normal form.

For instance the NTL term 
\[\Pi^{(1,1,2)} \Upsilon^{\{1,4\}\{2\}\{3\}} \Pi^{(1,3)} \Pi^{(1,1)} C D D D D C B C\]

reduces to the normal term

\[\Upsilon^{\{1,2\}\{3\}\{4\}\{5\}} \Pi^{(1,4)} C \Pi^{(1,1,1,1,1)} D \Pi^{(1,1,1,1,1)} D C I I I I I I I I \Pi^{(4,1,1,1,1)} D \Upsilon^{\{4\}\{1\}\{2\}\{3\}\{5\}\{6\}} \Pi^{(1,2,1,1,1)} D B C C I I I I I I\]

The following file is an example of how the tool is used on the above term to reach the normal term above using several different reduction sequences.\\

\url{https://github.com/owl77/natlog/blob/main/ntltest3.txt}\\

\end{example}

\section{Conclusion}

In this paper we have presented a formal system (NTL) based on fundamental aspects of the structure of natural language including notably the ability to express intensionality and the feature that quantification can be expressed without employing variable binding.  NTL comes with a system of reductions which can be interpreted as expressing how certain fundamental syntactical transformations preserve a basic logical aspect of meaning. This basic level of logical meaning can be uniquely determined and expressed via normalization (or equivalently via the corresponding Bealer term) thanks to Theorem 4.23.  We are in the presence of a purely syntactical characterization of the core logical meaning of a term. The NTL terms themselves can be interpreted as conveying finer grained meaning beyond such a core logical meaning, in the style of Bealer's  intentional logic T2\cite{QC}[p.64], though of course NTL is much more fine-grained than BL. It can be hoped that NTL may provide a useful tool not only for the formalization and analysis of natural language but also for elucidating classical philosophical problems in the spirit of Bealer's approach using his system T2 or extensions thereof\cite{QC}. An interesting application of NTL with its versality as an intensional logic could be to the study of how historical natural language based systems of logic expressed or might have expressed multiple generality, for instance  the logic of Stoicism \cite{mult} or the tradition of Boethius \cite{pin}.
In our presentation of NTL we employed for simplicity logical operators corresponding only to standard classical connectives and the universal quantifier but it is clear that the results of this paper can be directly generalized to any kinds of connective, modality and quantifier-like operators. We also focused on syntax and logical meaning preserving syntactical transformations  and have left out considerations on inference. Future work, in connection to the philosophical questions in \cite{QC}, will involve formulating and adding axioms and  rules of inference to NTL and investigating the resulting properties. Another avenue for future exploration could involve developing variants of NTL and the transformations $\beta$ and $\nu$ geared to second-order logic or extensions of second-order logic along the lines of the semantic version of a variant of a fragment of NTL used by Zalta in \cite{Zalta}.  Also it might be possible to adapt the results of this work to the $\lambda$-calculus to obtain a decomposition of terms of the simply typed lambda calculus analogous to $\beta$ but where the NTL constructors are now elements of the polymorphic system $F$. Since  the use of weaver theory in NTL  has some similarity to the use of permutations in operad theory it is conceivable that the framework of NTL may provide an extension of operad theory that includes generalized diagonalization. The author is also currently working on a graph-theoretic notation for NTL which involves a geometrical interpretation of the reductions. We sketch such a graph theoretic system in the appendix.

\section*{Acknowledgments}

The author received no financial support for the research contained in this paper.  The author has no conflict of interest to report.

\bibliographystyle{plain}
\bibliography{References}

\section*{Appendix: Sketch of a Graph Calculus for NTL}

We sketch here how the terms and reduction rules of Natural Term Logic can be given a geometric or topological visualization or interpretation. We also point out how such an interpretation can serve as a foundation for
a generalization of operad theory which we call combinatory operad theory.





\subsection*{GTC graphs}

As we saw,  NTL terms are built up from primitive terms and the application of constructors. 
We now proceed to describe what we call Geometric Term Calculus (GTC).  NTL terms correspond to GTC graphs.  Primitive terms of valency $n$ correspond to pictures as in fig. 1. A primitive term is thus represented by a isosceles triangle labelled with the symbol of the corresponding primitive term - which we call a \emph{basic triangle}. The top of the triangle is called the \emph{head} and the bottom lines (whose number corresponds to the valency of the primitive term) are called the \emph{legs}. 0-valent primitive terms correspond to triangles without legs. The exception for triangles with one leg is the unit $I$ which corresponds to a vertical line as in fig.1.

Every GTC graph can be considered (as we shall see in detail) as having a unique head and its set of legs. As such we also represent a general GTC graph by a labelled triangle in the same way we do for primitive terms with it being understood that such a triangle is an abstraction or a black box,  the contents of which are abstracted from or hidden (these are the grey triangles in fig.1). The label in this case is just a tag  and does not correspond to the label of a triangle of a primitive term.

GTC graphs are built up inductively as follows. Here we assume that the triangles represent arbitrary graphs.  Given a graph $C$ and a weaver $W$ we may apply  $W$ to the legs  of $C$ to obtain a graph as in fig 1. The grey triangles illustrate the we may again consider the resulting graph as an abstract triangle in which we ignore the internal content. If our graph $C$ corresponds to a NTL term $C$ then this construction corresponds to the NTL term $\Upsilon^W C$.

Given a graph $C$ with $n$ legs and $n$ other graphs (which can include the line $I$) we may plug the heads of these graphs into the legs of $C$ in different ways (according to the selector $S$) as exemplified in fig.1.  Plugging in $I$ should be understood as topologically doing nothing, merely extending the leg.
The crucial ingredient here is the blocking of a leg indicated by a small black circle. Once a leg has been blocked no graph can be plugged into it. The selector $S$ determines how many legs are blocked when we plug a certain graph into a certain leg (it is the number of legs of the graph to be plugged in minus the corresponding number of the selector). The blocking is always from the right to the left. Thus in the bottom construction in fig.1 the second leg receives graph $C$ with 3 legs and the second entry in the selector is $2$ so 1 leg of $C$ will be blocked when it is plugged in.

Fig.2 illustrates how the remaining NTL constructors are interpreted through circles.  In the corresponding 'plugging in' (which is not related to $\Pi$) there can never be any blocking of legs. The special triangle for $\Lambda_\forall$ has a peculiar form of plugging in in which the right-most leg is not blocked but is integrated backwards into the triangle itself.

We have now completed the description of how all GTC graphs are built up.  Fig.3 gives an example of how a generic GTC graph might look.

\subsection*{The topology of NTL reductions}

Using the correspondence developed in the previous section we can see how NTL reductions correspond to certain transformations applied to GTC graphs.  Interestingly several of the reductions appear to be topologically trivial or immediate at the level of GTC graphs (we consider the triangles and the ends of the legs fixed while the legs themselves can move in $n$-dimensional space without crossing each other). This appears to be the case for $R_{I}, C_\Upsilon, Id_\Upsilon, R_\Pi, P_\Pi$ and the other pushing-in reductions. See figs. 4a, 4b and 5. Particularly interesting is the case of $R_\Pi$ in which the reduction seems to be nothing more than 'yanking down' the blocked legs (fig. 5). \textcolor{red}{We do note use $R_{I2}$ is the version of NTL of the paper}.
From the point of view of GTC the most interesting reduction appears to be $R_\Upsilon$ which geometrically involves a kind of 'unravelling' of the legs of main graph of a $\Pi$-application (see fig. 6) in which argument graphs can be multiplied or 'cloned'.
The remaining pushing-in reductions for the logical constructors are either trivial or do not offer anything new from the GTC perspective.

\subsection*{Combinatory Operads}

For those familiar with the graphical notation for multicategories and operads it could be tempting to make a connection between GTC and the theory of operads.  GTC certainly suggests an interesting generalization of symmetric operads which we will describe below. However it seems highly unlikely that such a category theoretic setting could funish a faithful representation of NTL or GTC. The reason has to do with interpreting and capturing the blocking of legs in $\Pi$-applications.  Cartesian closed categories are the usual setting for operads and such categories correspond to the simply typed $\lambda$-calculus.  But NTL can represent Bealer Logic which is essentially untyped.  When we use $\Pi$-application and block legs we are really allowing differently typed terms in the same argument place, something which is impossible in the way multicomposition is defined in operad theory.  For instance we can compose $f: X^2 \rightarrow X$ with two copies of $g: X^2 \rightarrow X $ but not with one copy of $g$ and $g'(x) = \lambda y.g(x,y)$. The most we can do is to covertly indicate the $\lambda$-abstraction through something like the Henkin construction, replacing bound variables with special constants. In this case the blocking of a leg would correspond to composing with a special element $1 \rightarrow X$ which can be seen as an element of $P(0)$ for an operad $P$.

But let us now describe a generalization of symmetric operads. 

For simplicity we will only work over a cartesian closed category $\mathcal{C}$ (these considerations probably extend to symmetric monoidal categories). Let $X$ be an object in $\mathcal{C}$. As usual we consider, for $n\geq 1$, sets $P(n) = \{ f: X^n \rightarrow X\}$ (which in the multicategory perspective are represented as triangles similar to in GTC).  For a symmetric operad there is an action of the symmetric group $S_n$ on $P(n)$ as well as the usual properties for the unit and multi-composition. 
A combinatory operad is the same scenario but where we have now an action of $W(n,m)$ (weavers of order $(n,m)$) which takes a $P(n)$ to a $P(m)$.  This action is defined by generalized diagonalization using the cartesian closed structure of $\mathcal{C}$. For instance the action of $W = ((1,2,3), \{2,3\}\{1\})$ on a $f: X^3 \rightarrow X$ is defined as follows. Consider $\Delta: X \rightarrow X\times X$ defined by $\langle id_X,id_X\rangle$ and $\sigma : X\times(X\times X) \rightarrow (X\times X)\times X$ defined by $\langle \pi_{X\times X}, \pi_X \rangle$ and then the map $w = \sigma \circ (id_X \times \Delta): X^2 \rightarrow X^3$.  We then define the action of $W$ on $f$ as being $ f \circ w : X^2 \rightarrow  X$ which corresponds to twisting the legs of a GCT graph according to $\Upsilon^W$.  We can generalize this constrution for any weaver $W$.


 The NTL and GTC 'unit' is obviously the operad unit   $id_X$ in $P(1)$.

We might further propose that there is some way, for interpreting a general $\Pi$-application,  to assign to each blocked leg a distinguished morphism (tagged by the corresponding morphism) $1_{f}: 1 \rightarrow X$ which is appliced in the context of ordinary operad composition. The argument is certainly blocked in the sense that nothing more in a $P(n)$ can be attached to it. All this certainly demands further consideration.

Let us now interpret the NTL or GTC reductions in terms of combinatory operads. The reductions involving $I$ and some other are clearly just basic algebraic properties or closely related to such.

The most interesting reduction from the operad perspetive is again $R_\Upsilon$. Let us look at a simple case. Let $f: X^2 \rightarrow X$ be in $P(2)$ and consider the weaver $W = ((1,2), \{1,2\})$ which acts as ordinary diagonalization to yield $ W(f) = f \circ \Delta : X \rightarrow X\times X \rightarrow X$ in $P(1)$.
Now suppose we compose with $g : X^2 \rightarrow X$ to obtain $ W(f) \circ g  : X^2 \rightarrow X$ in $P(2)$. What does $R_\Upsilon$ say in this case? Consider the case in which $\mathcal{S}$ is the category of sets. Then $W(f) \circ g (x,y) =  f(g(x,y), g(x,y))$. This map can be obtained equivalently by first cconsidering
multicomposition of $f$ with two 'copies' of $g$ which can also be represented as $ f \circ (g \times g) : X^4 \rightarrow X$ in $P(4)$ and then applying the weaver $W = ((1,2,3,4), \{1,3\}\{2,4\})$ to this composition to obtain exactly the same map in $P(2)$. This very likely carries over to an CCC. Thus $R_\Upsilon$ expresses a very natural property to include in the abstract definition of a combinatory operad.\\


 \includegraphics[scale=0.4]{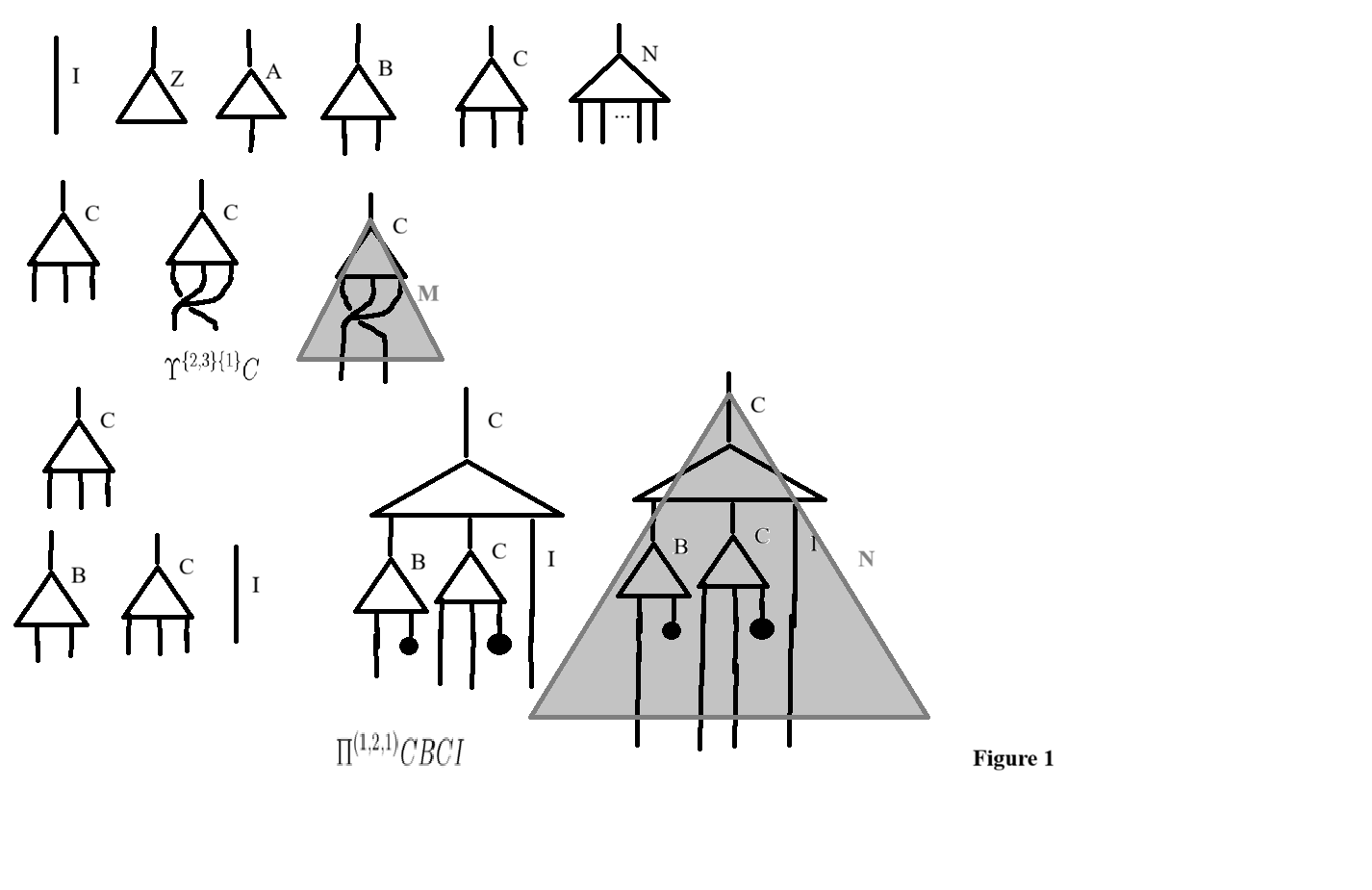}

 \includegraphics[scale=0.4]{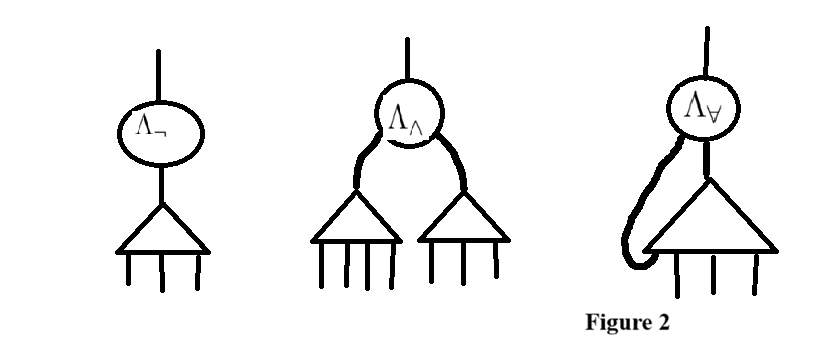}

 \includegraphics[scale=0.4]{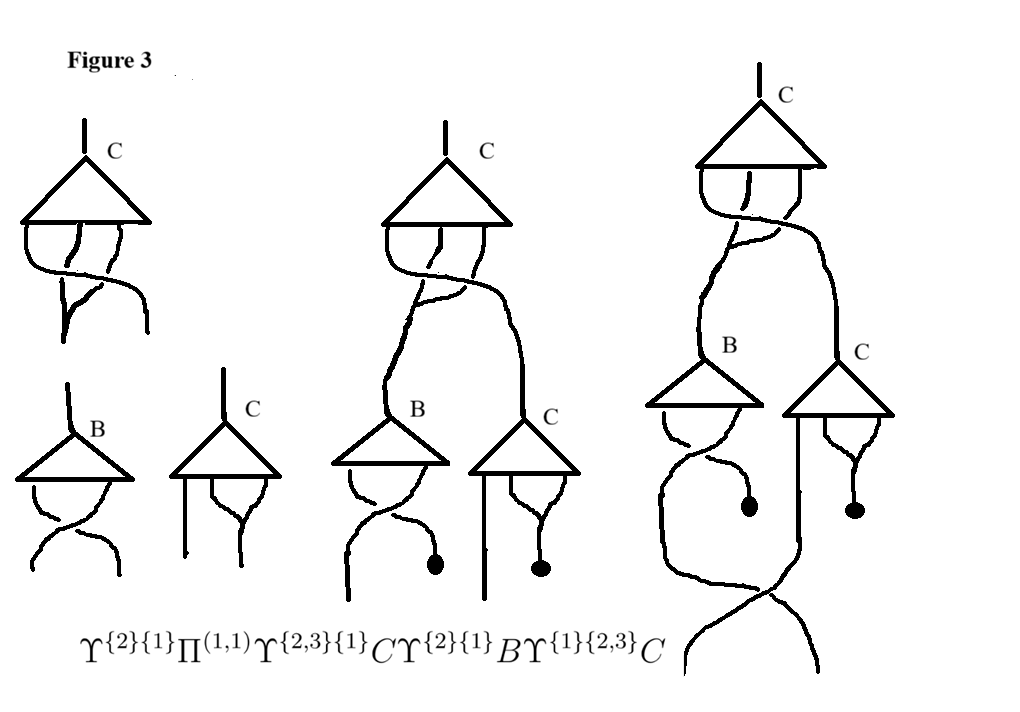}

 \includegraphics[scale=0.4]{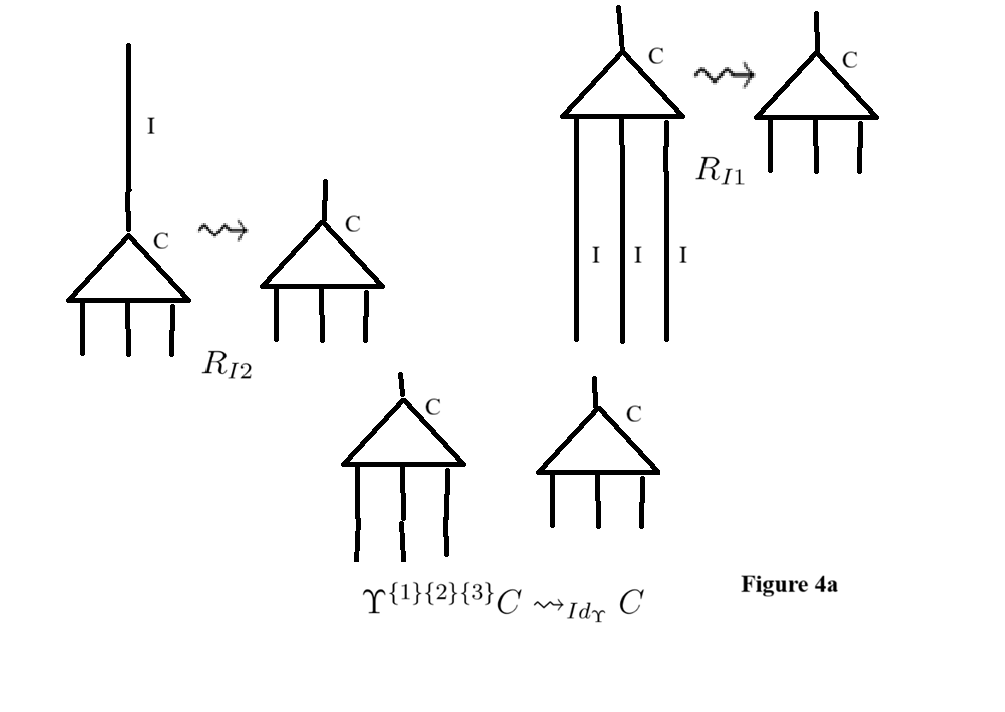}

 \includegraphics[scale=0.4]{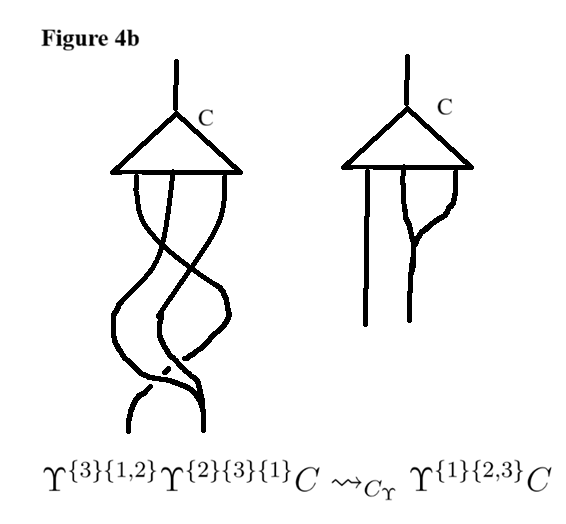}

 \includegraphics[scale=0.4]{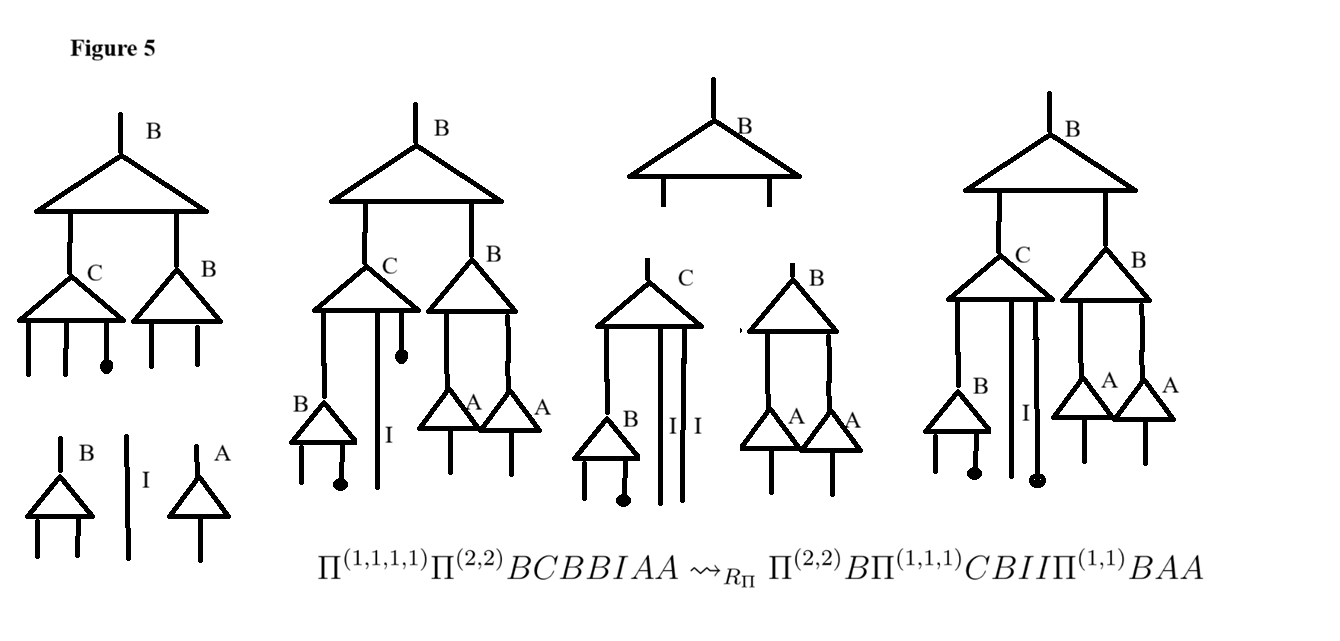}

 \includegraphics[scale=0.4]{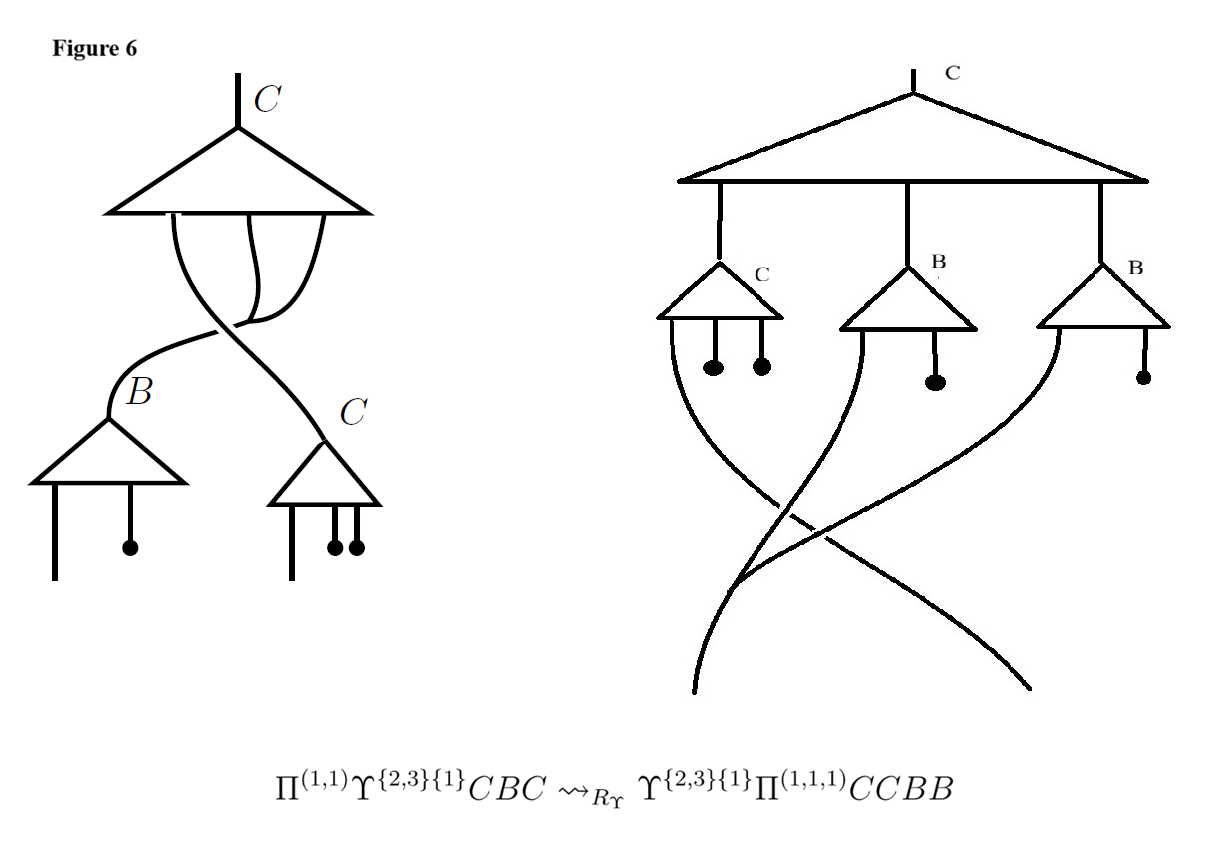}

\includegraphics[scale=0.4]{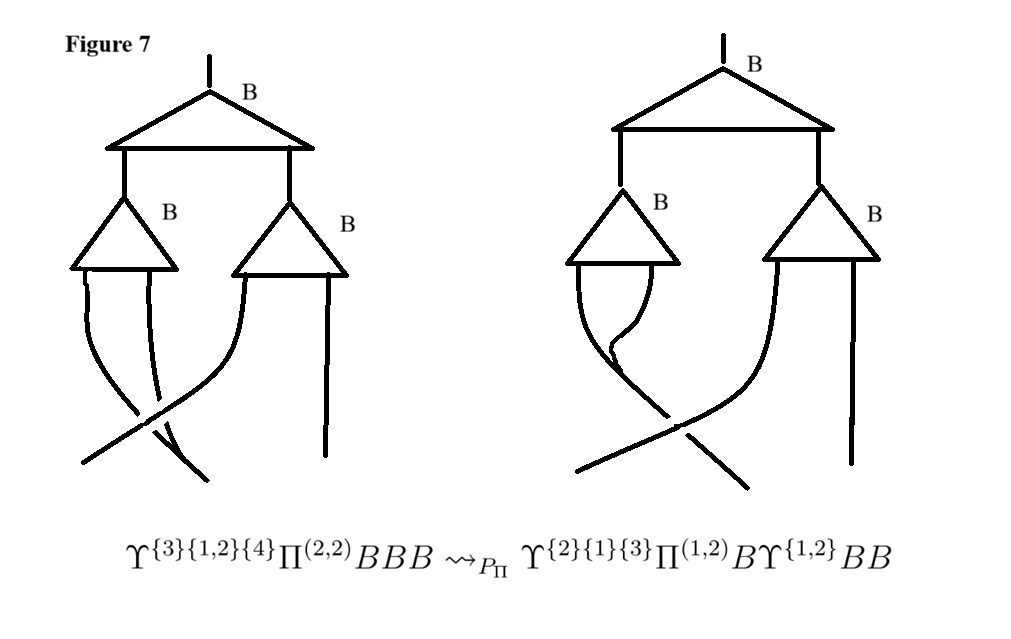}

\includegraphics[scale=0.4]{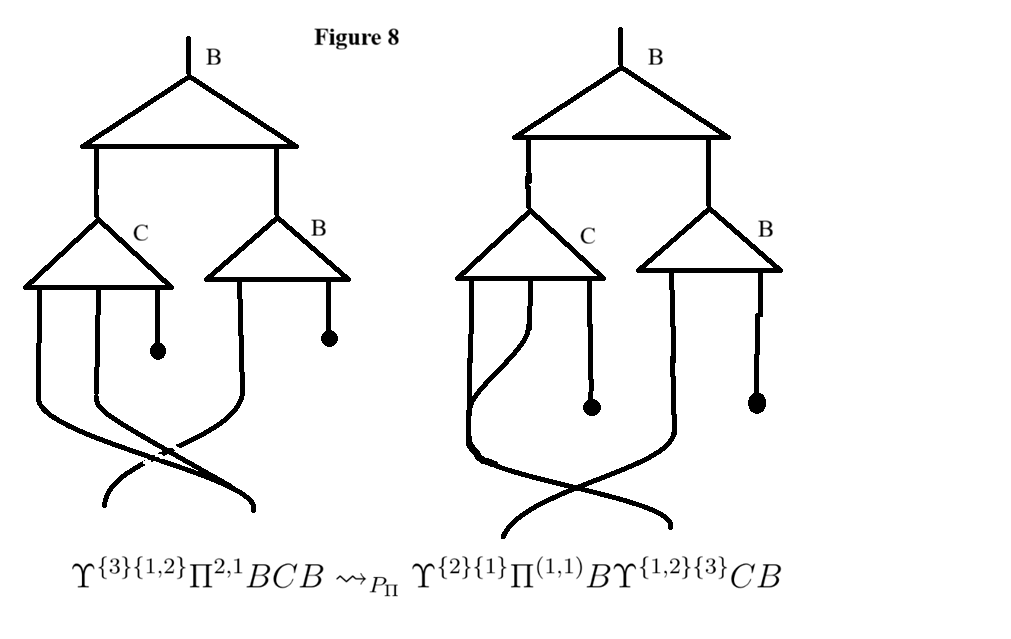}

\end{document}